\documentclass[12pt,reqno]{amsart}

\usepackage{amsfonts,color,amsthm,amsmath,amssymb}

\usepackage{color}
\allowdisplaybreaks[2]
\def\Box{\vcenter{\vbox{\hrule\hbox{\vrule
     \vbox to 8.8pt{\hbox to 10pt{}\vfill}\vrule}\hrule}}}

\def\qed{{\hfill$\square$}}
\def\proof{{\vspace{-0.0cm}\bf Proof: \,}}

\def\Z{{\mathbb Z}}
\def\Q{{\mathbb Q}}

\def\C{{\mathbb C}}
\def\F{{\mathbb F}}
\def\mod{{\mathrm{mod\,\,}}}

\def\Tr{{\mathrm{Tr}}}
\def\Norm{{\mathrm{Norm}}}

\def\ord{{\mathrm{ord}}}

\newtheorem{theorem}{Theorem}

\newtheorem{remark}[theorem]{Remark}
\newtheorem{corollary}[theorem]{Corollary}

\newtheorem{proposition}[theorem]{Proposition}
\newtheorem{example}[theorem]{Example}

\numberwithin{equation}{section}

\begin{document}
\title[Inequivalence of skew Hadamard difference sets]{Inequivalence of skew Hadamard difference sets and triple intersection numbers modulo a prime}
\author{Koji Momihara  }
\thanks{Mathematics Subject Classification (2010). Primary 05B10; Secondary 05E30.}
\thanks{K. Momihara is with the Faculty of Education, Kumamoto University,
2-40-1 Kurokami, Kumamoto 860-8555, Japan (e-mail:
momihara@educ.kumamoto-u.ac.jp). The work of
K. Momihara was supported by JSPS under Grant-in-Aid for Young Scientists (B) 25800093 and Scientific Research (C) 24540013.}
\keywords{skew Hadamard difference set, Feng-Xiang difference set, Paley difference set}

\maketitle

\begin{abstract}
Recently, Feng and Xiang \cite{FX113} found a new construction of skew Hadamard difference sets in elementary abelian groups. In this paper, we introduce a new invariant for equivalence of skew Hadamard difference sets, namely triple intersection numbers modulo a prime, and discuss inequivalence between Feng-Xiang skew Hadamard difference sets and the Paley difference sets. As a consequence, we show that their construction produces infinitely many  
skew Hadamard difference sets inequivalent to the Paley difference sets. 
\end{abstract}
\section{Introduction}
Let $(G,+)$ be an  (additively written) finite  group of order $v$. A $k$-subset $D\subseteq G$ is called a {\it $(v,k,\lambda)$ difference set} if the list of differences $x-y$ with  $x,y\in D$ and $x\not=y$ covers each nonzero element of $G$ exactly $\lambda$ times. We say that two difference sets $D_1$ and $D_2$ with the same parameters in an abelian group $G$ are {\it equivalent} if there exists an automorphism $\sigma\in \mathrm{Aut}(G)$ and an element $x\in G$ such that $\sigma(D_1)+x=D_2$. For general theory on difference sets, we refer the reader to \cite{Beth}.

A difference set $D$ in a finite group $G$ is called {\it skew Hadamard} if $G$ is the disjoint union of $D$, $-D$, and $\{0\}$. The primary example (and for many years, the only known example in abelian groups) of skew Hadamard difference sets is the classical Paley difference set in $(\F_q,+)$ consisting of the nonzero squares of $\F_q$, where $\F_q$ is the finite field of order $q$, and $q$ is a prime power congruent to $3$ modulo 4.  This situation changed dramatically in recent years. Skew Hadamard difference sets are currently under intensive study;  see  
\cite{CP,DAW,DY06,DWX07,F11,FMX11,FX113,M,WH09,WQWX} for recent results and the introduction of \cite{FX113} for a short survey of  skew Hadamard difference sets and related problems.
 
There were two major 
conjectures in this area: (i) If an abelian group $G$
contains a skew Hadamard difference set, then $G$ is necessarily elementary
abelian. (ii) Up to equivalence the Paley difference sets mentioned above
are the only skew Hadamard difference sets in abelian groups. The former
conjecture is still open in general. 
The latter conjecture turned out to be false: Ding
and Yuan \cite{DY06} constructed a family of skew Hadamard difference sets in
$(\F_{3^m},+)$, where $m\ge 3$ is odd, by using Dickson polynomials of order $5$ and showed that two examples
in the family are inequivalent to the Paley difference sets. Very recently, Ding, Pott, and Wang \cite{DAW} found more skew Hadamard difference sets inequivalent to the Paley difference sets from 
Dickson polynomials of order $7$. 
Muzychuk \cite{M} gave a prolific construction of skew Hadamard difference sets in an elementary abelian group of order $q^3$  and showed that  his skew Hadamard
difference sets are inequivalent to the Paley difference sets.
Although many other constructions have been known recently, as far as the author knows, there has been no theoretical result on the inequivalence problem of  skew Hadamard difference sets except for \cite{M}. Indeed, in most of recent papers, the inequivalence of skew Hadamard difference sets were checked by computer. 
Here, we should remark that some of known invariants for equivalence of ordinary difference sets (e.g., $p$-ranks of the symmetric designs developed from difference sets) do not contribute anything to  the inequivalence problem of skew Hadamard difference sets $D$ since they are determined from only  the parameters of $D$ \cite[Chapter VI, Theorem 8.22]{Beth}. 
On the other hand, some invariants (e.g., triple intersection numbers)
are difficult to compute without computer. 

A classical method for constructing  difference sets in the additive groups of finite fields is to use cyclotomic classes of finite fields. 
Let $p$ be a prime,  $f$ a positive integer, and let $q=p^f$. Let $k>1$ be an integer such that $N|(q-1)$, and $\gamma$ be a primitive element of $\F_q$. 
Then the cosets $C_i^{(N,q)}=\gamma^i \langle \gamma^N\rangle$, $0\leq i\leq N-1$, are called the {\it cyclotomic classes of order $N$} of $\F_q$. 

In this paper, we are particularly interested in the following construction of 
skew Hadamard difference sets given by Feng and Xiang~\cite{FX113}. 
\begin{theorem}\label{FXtheorem}{\em (\cite[Theorem 3.2]{FX113})}
Let $p_1\equiv 7\,(\mod{8})$ be a prime, $N=2p_1^m$, and 
let $p\equiv 3\,(\mod{4})$ be a prime such that $f:=\ord_N(p)=\phi(N)/2$. Let $s$ be any odd integer, $I$ any subset of $\Z/N\Z$ such that 
$\{i\,(\mod{p_1^m})\,|\,i\in I\}=\Z/p_1^m\Z$, and  let $D=\bigcup_{i\in I}C_i^{(N,q)}$, where $q=p^{fs}$. Then, $D$ is a skew Hadamard 
difference set. 
\end{theorem}
We call the difference set in Theorem~\ref{FXtheorem} as the {\it Feng-Xiang skew Hadamard difference set} with index set $I$.  Let $t$ be an odd integer and 
$\gamma$ be a primitive element of $\F_{q^t}$. Put $\omega=
\gamma^{(q^t-1)/(q-1)}$.
If $D=\bigcup_{i\in I}C_i^{(N,q)}=\bigcup_{i\in I}\omega^i \langle \omega^N\rangle$ is a Feng-Xiang skew Hadamard difference set, then so does 
$D'=\bigcup_{i\in I}C_i^{(N,q^t)}=\bigcup_{i\in I}\gamma^i\langle \gamma^N\rangle$. We call $D'$ the {\it lift} of $D$ to $\F_{q^t}$. Furthermore, throughout this paper,  we denote the set $\bigcup_{i\in I}\omega^{ti} \langle \omega^N\rangle$ by $D^{(t)}$. It is clear that 
 $D^{(t)}$ is also a Feng-Xiang skew Hadamard difference set if $\gcd{(t,N)}=1$. 

In this paper, we introduce a new invariant for equivalence of skew Hadamard difference sets, namely triple intersection numbers modulo a prime, and show that infinitely many Feng-Xiang skew Hadamard difference sets are inequivalent to the
Paley difference sets by using ``recursive'' techniques.  Besides 
the existence of infinitely many skew Hadamard difference sets  inequivalent 
to the Paley difference sets, our technique may contribute to inequivalence problems  on combinatorial objects defined in finite fields not only on skew Hadamard difference sets. 

This paper is organized as follows. In Section~\ref{Sec2}, we introduce some 
preliminaries on characters of finite fields, and present a proposition on  divisibility of a character sum over a finite field by its characteristic. In Section~\ref{Sec3}, we introduce the concept of  ``triple intersection numbers modulo a prime'' and we give two sufficient conditions for 
lifts of  Feng-Xiang 
skew Hadamard difference sets being inequivalent to the Paley difference sets. As an example, we show that there are infinitely many integers $t$ such that the lifts of a Feng-Xiang skew Hadamard difference set in $\F_{11^{3}}$ to $\F_{11^{3t}}$ are inequivalent to the Paley difference sets. In Section~\ref{Sec4}, we conclude this paper with further examples of  
skew Hadamard difference sets inequivalent  to the Paley difference sets, and some open problems. 
\section{Preliminaries on characters}\label{Sec2}
Let $p$ be a prime, $f$ a positive integer, and $q=p^f$. The canonical additive character $\psi$ of $\F_q$ is defined by 
$$\psi\colon\F_q\to \C^{\ast},\qquad\psi(x)=\zeta_p^{\Tr _{q/p}(x)},$$
where $\zeta_p={\rm exp}(\frac {2\pi i}{p})$ and $\Tr _{q/p}$ is the trace from $\F_q$ to $\F_p$. For a multiplicative character 
$\chi_N$ of order $N$ of $\F_q$, we define the {\it Gauss sum} 
\[
G_q(\chi_N)=\sum_{x\in \F_q^\ast}\chi_N(x)\psi(x), 
\] 
which belongs to the ring $\Z[\zeta_{pN}]$ of integers in the cyclotomic field $\Q(\zeta_{pN})$.
Let $\sigma_{a,b}$ be the automorphism of $\Q(\zeta_{pN})$ determined 
by 
\[
\sigma_{a,b}(\zeta_N)=\zeta_{N}^a, \qquad
\sigma_{a,b}(\zeta_p)=\zeta_{p}^b 
\]
for $\gcd{(a,N)}=\gcd{(b,p)}=1$. 
Below are several basic properties of Gauss sums \cite{BEW97}: 
\vspace{0.2cm}
\begin{itemize}
\item[(i)] $G_q(\chi_N)\overline{G_q(\chi_N)}=q$ if $\chi$ is nontrivial;
\item[(ii)] $G_q(\chi_N^p)=G_q(\chi_N)$, where $p$ is the characteristic of $\F_q$; 
\item[(iii)] $G_q(\chi_N^{-1})=\chi_N(-1)\overline{G_q(\chi_N)}$;
\item[(iv)] $G_q(\chi_N)=-1$ if $\chi_N$ is trivial;
\item[(v)] $\sigma_{a,b}(G_q(\chi_N))=\chi_N^{-a}(b)G_q(\chi_N^a)$.
\end{itemize}
\vspace{0.2cm}
In general, the explicit evaluation of  Gauss sums is a very difficult problem. There are only a few cases where the Gauss sums have been evaluated. 
The most well known case is {\it quadratic} case, i.e., the case where $N=2$. 
In this case,  it holds that 
\begin{equation}\label{eq:quad}
G_f(\chi_N)=(-1)^{f-1}\bigg(\sqrt{(-1)^{\frac{p-1}{2}}p}\bigg)^f, 
\end{equation}
cf.~\cite[Theorem 11.5.4]{BEW97}. 
The next simple case is the so-called {\it semi-primitive case} (also 
referred to as {\it uniform cyclotomy} or {\it pure Gauss sum}), where there 
exists an integer $j$ such that $p^j\equiv -1\,(\mod{N})$, where $N$ is the order of 
the multiplicative character involved. The explicit evaluation of Gauss sums in this case is given in \cite{BEW97}. 
The next interesting case is the {\it index $2$} case, where the subgroup $\langle p\rangle$ generated by $p\in (\Z/{N}\Z)^\ast$ is of index $2$ in $(\Z/{N}\Z)^\ast$ and $-1\not\in \langle p\rangle $. 
A complete solution to the problem of evaluating index $2$ Gauss sums 
was recently given in \cite{YX10}. The following is the result on evaluation of index $2$ Gauss sums, which 
was used to prove Theorem~\ref{FXtheorem}. 
\begin{theorem}\label{Sec2Thm2}{\em (\cite{YX10}, Case D; Theorem~4.12)}
Let $N=2p_1^m$, where $p_1>3$ is a prime such that $p_1\equiv 3\,(\mod{4})$ and $m$ is a positive integer. Assume that $p$ is a 
prime such that  $[(\Z/N\Z)^\ast:\langle p \rangle]=2$. 
Let $f=\phi(N)/2$, $q=p^f$, and $\chi$ be a multiplicative character of order $N$ of $\F_q$. Then, for 
$0\le t\le m-1$, we have
\begin{eqnarray*}
G_q(\chi^{p_1^t})&=&\left\{
\begin{array}{ll}
(-1)^{\frac{p-1}{2}(\frac{f-1}{2}-1)}p^{\frac{f-1}{2}-hp_1^t}
\sqrt{p^\ast}\bigg(\frac{b+c\sqrt{-p_1}}{2}\bigg)^{2p_1^t},&  \mbox{if $p_1\equiv 3\,(\mod{8})$,}\\
(-1)^{\frac{p-1}{2}\frac{f-1}{2}}p^{\frac{f-1}{2}}
\sqrt{p^\ast},&  \mbox{if $p_1\equiv 7\,(\mod{8})$;}
 \end{array}
\right.\\
G_q(\chi^{2p_1^t})&=&p^{\frac{f-p_1^t h}{2}}\bigg(\frac{b+c\sqrt{-p_1}}{2}\bigg)^{p_1^t};\\
G_q(\chi^{p_1^m})&=&(-1)^{\frac{p-1}{2}\frac{f-1}{2}}p^{\frac{f-1}{2}}\sqrt{p^\ast}, 
\end{eqnarray*}
where $p^\ast=(-1)^{\frac{p-1}{2}}p$, $h$ is the class number of $\Q(\sqrt{-p_1})$, and $b$ and $c$ are integers 
determined by $4p^{h}=b^2+p_1c^2$ and $bp^{\frac{f-h}{2}}\equiv -2\,(\mod{p_1})$. 
\end{theorem}
The following theorem, called the {\it Davenport-Hasse lifting formula}, is  useful for evaluating Gauss sums. 
\begin{theorem}\label{thm:lift}{\em (\cite[Theorem~5.14]{LN97})}
Let $\chi$ be a nontrivial multiplicative character of $\F_q=\F_{p^f}$ and
let $\chi'$ be the lifted character of $\chi$ to the extension field $\F_{q'}=\F_{p^{fs}}$, i.e., $\chi'(\alpha):=\chi(\Norm_{\F_{q'}/\F_q}(\alpha))$ for any $\alpha\in \F_{q'}^\ast$.
It holds that
\[
G_{q'}(\chi')=(-1)^{s-1}(G_q(\chi))^s.
\]
\end{theorem}
In relation to Gauss sums,  we need to define the Jacobi sums. 
Let $\chi$ and $\chi'$ be multiplicative characters of $\F_q$. We  define the sum 
\[
J(\chi,\chi')=\sum_{x\in \F_q,x\not=0,1}\chi(x)\chi'(1-x),
\]
the so-called {\it Jacobi sum} of $\F_q$. It is known \cite[Theorem 5.21]{LN97} that 
if $\chi$, $\chi'$, and $\chi\chi'$ are nontrivial, it  holds that 
\begin{equation}\label{GaussJacobi}
J(\chi,\chi')=\frac{G_{q}(\chi)G_{q}(\chi')}{G_{q}(\chi\chi')}. 
\end{equation}
We will use this formula later. 

Now we are interested in computing the following character sum modulo the characteristic $p$: 
\[
\sum_{x\in \F_{p^f}}\chi(x^{i_1}(x+1)^{i_2}(x+a)^{i_3}) \, (\mod{p}), 
\]
where $\chi$ is a multiplicative character of $\F_{p^f}$. 
The following theorem is well known as the {\it Weil theorem on multiplicative 
character sums}. 
\begin{theorem}\label{multiweil}{\em (\cite[Theorems 5.39 and 5.41]{LN97})}
Let $\chi$ be a multiplicative character of $\F_q$ of order $N>1$ and 
$f\in \F_q[x]$ be a monic polynomial of positive degree that is not an 
$N$th power of a polynomial. Let $d$ be the number of distinct roots of $f$ 
in its splitting field over $\F_q$ and suppose that $d\ge 2$. Then there 
exist complex numbers $w_1,\ldots,w_{d-1}$, only depending on $f$ and $\chi$, such that for any positive integer $t$ we have 
\begin{equation}\label{weil1}
\sum_{x\in \F_{q^t}}\chi'(f(x))=-w_1^t-\cdots-w_{d-1}^t,
\end{equation}
where $\chi'$ is the lift of $\chi$ to $\F_{q^t}$. In particular, it holds that 
\begin{equation}\label{weil2}
\bigg|\sum_{x\in \F_{q^t}}\chi'(f(x))\bigg|\le (d-1)\sqrt{q^t}.  
\end{equation}
\end{theorem}
\vspace{0.5cm}
With the notations above, we set $d=3$. By Warning's formula~\cite[Theorem 1.76]{LN97}, $w_1^t+w_2^t$ can be expressed 
as follows:
\begin{equation}\label{war1}
w_1^t+w_2^t=\sum_{j=0}^{\lfloor t/2\rfloor}(-1)^j\frac{t}{t-j}{t-j \choose j}(w_1+w_2)^{t-2j}(w_1w_2)^j,  
\end{equation}
where each coefficient of $(w_1+w_2)^{t-2j}(w_1w_2)^j$ takes an integral value. 
We assume that $f(x)$ can be decomposed as $f(x)=(x-a)^{i_1}(x-b)^{i_2}(x-c)^{i_3}\in \F_q[x]$ for distinct $a,b,c\in \F_q$ and $i_1,i_2,i_3\not\equiv 0\,(\mod{N})$. 
As found in the proof of Theorem~5.39 (also Eqs. (5.19) and (5.23)) in \cite{LN97}, we have 
\begin{equation}\label{polydeg2}
w_1w_2=\sum_{g\in \Phi_2}\lambda(g), 
\end{equation}
where $\Phi_2$ is the set of all monic polynomials of degree $2$ over $\F_q$ and 
$\lambda$ is defined by $\lambda(g)=
\chi^{i_1}(g(a))\chi^{i_2}(g(b))\chi^{i_3}(g(c))$ for $g(x)\in \Phi_2$. 
Note that if $t$ is an odd prime, by Eq.~(\ref{war1}) we obviously have 
\[
w_1^t+w_2^t\equiv (w_1+w_2)^{t}\, (\mod{t}) 
\] 
since $w_1+w_2,w_1w_2\in \Z[\zeta_N]$. 

We will use the following proposition in the next section. 
\begin{proposition}\label{p-reduc}
Let $\chi$ be a multiplicative character of $\F_q$ of order $N>1$ and 
$f(x)=x^{i_1}(x+1)^{i_2}(x+a)^{i_3}\in \F_q[x]$ for $a\in \F_p^\ast\setminus \{1\}$ and  $i_1,i_2,i_3\not\equiv 0\,(\mod{N})$.  Assume that $p$ divides $J(\chi^{i_2},\chi^{i_3})J(\chi^{i_1},\chi^{i_2i_3})$ for $i_2$ and $i_3$ such that  $\chi^{i_2i_3}$ is nontrivial. Then, for any odd integer $t$ it holds that
\[
\sum_{x\in \F_{q^t}}\chi'(f(x))\equiv \bigg(\sum_{x\in \F_q}\chi(f(x))\bigg)^t \,(\mod{p}),
\]  
where $\chi'$ is the lift of $\chi$ to $\F_{q^t}$. 
\end{proposition}
\proof
By Eq.~(\ref{weil1}) of Theorem~\ref{multiweil}  and Eq.~(\ref{war1}), it is enough to show that 
$w_1w_2\equiv 0\,(\mod{p})$. By Eq.~(\ref{polydeg2}), we need to compute 
the sum $\sum_{g\in \Phi_2}\lambda(g)$. By the definition of $\lambda$, we have 
\begin{subequations}
\begin{align}
\sum_{g\in \Phi_2}\lambda(g)&=\sum_{x\in \F_q}\sum_{y\in \F_q}
\chi^{i_1}(y)\chi^{i_2}(1-x+y)\chi^{i_3}(a^2-ax+y)\nonumber\\
&=\sum_{z\in \F_q}\sum_{y\in \F_q}
\chi^{i_1}(y)\chi^{i_2}(z)\chi^{i_3}(a^2+y(1-a)-a+az)\nonumber\\
&=\sum_{z\in \F_q}\sum_{y\in \F_q\setminus \{a\}}
\chi^{i_1}(y)\chi^{i_2}(z)\chi^{i_3}(a^2+y(1-a)-a)\chi^{i_3}\bigg(1+\frac{a}{a^2+y(1-a)-a}z\bigg)\label{func1}\\
&\hspace{1cm}+\sum_{z\in \F_q}
\chi^{i_1}(a)\chi^{i_2}(z)\chi^{i_3}(az). \label{func2}
\end{align}
\end{subequations}
It is clear that 
\[
(\ref{func2})=\left\{ 
\begin{array}{ll}
\chi^{i_1}(a)\chi^{-i_2}(a)(q-1)&\quad  \text{if $\chi^{i_2i_3}$ is trivial,} \\
0 &\quad \text{if $\chi^{i_2i_3}$ is nontrivial.}
\end{array}
\right. 
\]
Next, we compute the sum (\ref{func1}). 
By the definition of Jacobi sums, we have 
\begin{align}
(\ref{func1})&=J(\chi^{i_2},\chi^{i_3})\sum_{y\in \F_q\setminus \{a\}}
\chi^{i_1}(y)\chi^{-i_2}\bigg(\frac{-a}{a^2+y(1-a)-a}\bigg)\chi^{i_3}(a^2+y(1-a)-a)\nonumber\\
&=J(\chi^{i_2},\chi^{i_3})\chi^{-i_2}(-a)\sum_{y\in \F_q\setminus \{a\}}
\chi^{i_1}(y)\chi^{i_2i_3}(a^2+y(1-a)-a)\nonumber\\
&=J(\chi^{i_2},\chi^{i_3})\chi^{-i_2}(-a)\chi^{i_2i_3}(a^2-a)\chi^{-i_1}\bigg(\frac{a-1}{a^2-a}\bigg)\sum_{w\in \F_q\setminus \{-1\}}
\chi^{i_1}(-w)\chi^{i_2i_3}(1+w). \label{jacobitra}
\end{align}
If $\chi^{i_2}\chi^{i_3}$ is nontrivial, we have
\[
(\ref{jacobitra})=\chi^{-i_2}(-a)\chi^{i_2i_3}(a^2-a)\chi^{i_1}(a)J(\chi^{i_2},\chi^{i_3})J(\chi^{i_1},\chi^{i_2i_3}). 
\]
If $\chi^{i_2}\chi^{i_3}$ is trivial, by $J(\chi^{i_2},\chi^{i_3})=-\chi^{i_2}(-1)$, we have 
\begin{align*}
(\ref{jacobitra})&=\chi^{-i_2}(-a)\chi^{i_1}(a)J(\chi^{i_2},\chi^{i_3})\sum_{w\in \F_q\setminus \{-1\}}
\chi^{i_1}(-w)\\
&=-\chi^{-i_2}(-a)\chi^{i_1}(a)J(\chi^{i_2},\chi^{i_3})=\chi^{-i_2}(a)\chi^{i_1}(a). 
\end{align*}
By the assumption that  $p\,|\,J(\chi^{i_2},\chi^{i_3})J(\chi^{i_1},\chi^{i_2i_3})$ for $i_2$ and $i_3$ such that  $\chi^{i_2i_3}$ is nontrivial, 
we finally obtain
\begin{align*}
\sum_{g\in \Phi_2}\lambda(g)&=\left\{ 
\begin{array}{ll}
\chi^{i_1}(a)\chi^{-i_2}(a)q&\quad  \text{if $\chi^{i_2i_3}$ is trivial,} \\
\chi^{-i_2}(-a)\chi^{i_2i_3}(a^2-a)\chi^{i_1}(a)J(\chi^{i_2},\chi^{i_3})J(\chi^{i_1},\chi^{i_2i_3})&\quad \text{if $\chi^{i_2i_3}$ is nontrivial,}
\end{array}
\right. \\
&\equiv 0\, \pmod{p}. 
\end{align*}
Then the proof is complete.  \qed
\begin{remark}
{\em Let $q$ and $N$ be defined as in Theorem~\ref{FXtheorem}. Assume that  $i_1,i_2,$ and $i_3$ are odd and $\chi^{i_2i_3}$ is nontrivial. Note that $\chi^{i_1i_2i_3}$ is nontrivial since $i_1,i_2,$ and $i_3$ are odd. Then, by Eq.~(\ref{GaussJacobi}), we have 
\[
J(\chi^{i_2},\chi^{i_3})J(\chi^{i_1},\chi^{i_2i_3})=\frac{G(\chi^{i_2})G(\chi^{i_3})}{G(\chi^{i_2i_3})}\cdot \frac{G(\chi^{i_1})G(\chi^{i_2i_3})}{G(\chi^{i_1i_2i_3})}=\frac{G(\chi^{i_1})G(\chi^{i_2})G(\chi^{i_3})}{G(\chi^{i_1i_2i_3})}. 
\]
By Theorems~\ref{Sec2Thm2} and ~\ref{thm:lift}, we have 
\[
J(\chi^{i_2},\chi^{i_3})J(\chi^{i_1},\chi^{i_2i_3})=\bigg((-1)^{\frac{p-1}{2}\frac{f-1}{2}}p^{\frac{f-1}{2}}
\sqrt{p^\ast}\bigg)^{2s}\equiv 0\,(\mod{p}),  
\]
i.e., the condition of Proposition~\ref{p-reduc} is satisfied.  
}\end{remark}
\section{Triple intersection numbers modulo a prime} \label{Sec3}
Let  $D\subseteq \F_{q}$ be a skew Hadamard difference set and $\omega$ a primitive element of $\F_q$.  
For $a\in \F_p^\ast\setminus\{1\}$, let
\[
T_{\omega^\ell,a}(D):=|D\cap (D-\omega^\ell) \cap (D-a\omega^\ell)|. 
\]
The set $\{T_{\omega^\ell,a}(D)\,|\,0\le \ell\le q-2\}$ is an invariant for  equivalence of skew Hadamard difference sets. In fact, if $D'$ is a skew Hadamard difference set equivalent to $D$, namely $\sigma(D)=D'+x$ for $\sigma\in \mathrm{Aut}(\F_q,+)$ and $x\in \F_q$, we have 
\begin{align*}
\{T_{\omega^\ell,a}(D)\,|\,0\le \ell\le q-2\}
&=\{|\sigma(D\cap (D-\omega^\ell) \cap (D-a\omega^\ell))|:0\le \ell\le q-2\}\\
&=\{|\sigma(D)\cap (\sigma(D)-\sigma(\omega^\ell)) \cap (\sigma(D)-a\sigma(\omega^\ell))|:0\le \ell\le q-2\}\\
&=\{|D'\cap (D'-\omega^\ell) \cap (D'-a\omega^\ell)|:0\le \ell\le q-2\}. 
\end{align*}

If $D$ is the Paley difference set in $\F_q$, then $|\{T_{\omega^\ell,a}(D)\,|\,0\le \ell\le q-2\}|\le 2$ since 
$T_{1,a}(D)=T_{\omega^{2\ell},a}(D)$ and  $T_{\omega,a}(D)=T_{\omega^{2\ell+1},a}(D)$ 
for all $0\le \ell\le (q-3)/2$. Hence, if a skew Hadamard difference set $D'$ satisfies $|\{T_{\omega^\ell,a}(D')\,|\,0\le \ell\le q-2\}|\ge 3$, then 
$D'$ is inequivalent to the Paley difference set $D$.

Let  $D=\bigcup_{i\in I}C_i^{(N,q)}\subseteq \F_{p^{fs}}=\F_q$ be a Feng-Xiang skew Hadamard difference set.  It is clear that 
\[
|\{T_{\omega^\ell,a}(D)\,|\,0\le \ell\le q-2\}|=|\{T_{\omega^\ell,a}(D)\,|\,0\le \ell\le N-1\}|. 
\]
In this section, we compute the size of the set 
$
\{T_{\omega^\ell,a}(D)\, (\mod{t})\,|\,0\le \ell\le N-1\}
$
for a prime $t$. It is clear that this set is also an invariant for equivalence of skew Hadamard difference sets. Hence, if the set above contains at least three numbers, then 
$D$ is inequivalent to the Paley difference set.

Let $\chi_N$ be the multiplicative character of order $N$ of $\F_q$ such that $\chi_N(\omega)=\zeta_N$ and let $\eta_p$ be the quadratic character of $\F_p$. Note that  $\chi_N|_{\F_p}=\eta_p$. 
Since the characteristic function of $D=\bigcup_{i\in I}C_i^{(N,q)}$ is given by 
\[
f(x)=\frac{1}{N}\sum_{h\in I}\sum_{i=0}^{N-1}\zeta_N^{-i h}\chi_N^{i}(x), 
\]
we have 
\begin{align*}
&\hspace{0.6cm} N^3\cdot  T_{\omega^\ell,a}(D)\\
&=\, \sum_{x\in \F_{q}\setminus\{0,-1,-a\}}\sum_{h_1,h_2,h_3\in I}\bigg(\sum_{i_1=0}^{N-1}\zeta_N^{-i_1 h_1}\chi_N^{i_1}(x)\bigg)\bigg(\sum_{i_2=0}^{N-1}\zeta_N^{-i_2 h_2}\chi_N^{i_2}(x+\omega^\ell)\bigg)\bigg(\sum_{i_3=0}^{N-1}\zeta_N^{-i_3 h_3}\chi_N^{i_3}(x+a\omega^\ell)\bigg)\\
&=\, \sum_{x\in \F_{q}\setminus\{0,-1,-a\}}\sum_{h_1,h_2,h_3\in I}\sum_{i_1,i_2,i_3=0}^{N-1}\zeta_N^{-i_1 h_1-i_2 h_2-i_3 h_3}\chi_N^{i_1}(x)\chi_N^{i_2}(x+\omega^\ell)\chi_N^{i_3}(x+a\omega^\ell).
\end{align*}
Write $M=N/2$. 
By noting that  $\sum_{h\in I}\zeta_{N}^{2h}=0$, the above is expanded as follows:  
{\footnotesize \begin{align*}
&M\sum_{h_2,h_3\in I}\sum_{i_2,i_3=0}^{N-1}\zeta_N^{-i_2 h_2-i_3 h_3+\ell(i_2+i_3)}\sum_{x\in \F_{q}\setminus\{-1,-a\}}\chi_N^{i_2}(x+1)\chi_N^{i_3}(x+a)-M\sum_{h_2,h_3\in I}\sum_{i_2,i_3\in A}\zeta_N^{-i_2 h_2-i_3 h_3+\ell(i_2+i_3)}\eta_p(a)\\
&+M\sum_{h_1,h_3\in I}\sum_{i_1,i_3=0}^{N-1}\zeta_N^{-i_1 h_1-i_3 h_3+\ell(i_1+i_3)}\sum_{x\in \F_{q}\setminus\{0,-a\}}\chi_N^{i_1}(x)\chi_N^{i_3}(x+a)
-M\sum_{h_1,h_3\in I}\sum_{i_1,i_3\in A}\zeta_N^{-i_1 h_1-i_3 h_3+\ell(i_1+i_3)}\eta_p(-a+1)\\
&+M\sum_{h_1,h_2\in I}\sum_{i_1,i_2=0}^{N-1}\zeta_N^{-i_1 h_1-i_2 h_2+\ell(i_1+i_2)}\sum_{x\in \F_{q}\setminus\{0,-1\}}\chi_N^{i_1}(x)\chi_N^{i_2}(x+1)
-M\sum_{h_1,h_2\in I}\sum_{i_1,i_2\in A}\zeta_N^{-i_1 h_1-i_2 h_2+\ell(i_1+i_2)}\eta_p(a^2-a)
\\
&-M^2\sum_{h_1\in I}\sum_{i_1=0}^{N-1}\zeta_N^{-i_1 h_1+\ell i_1}\sum_{x\in \F_{q}\setminus\{0\}}\chi_N^{i_1}(x)
+M^2\sum_{h_1\in I}\sum_{i_1\in A}\zeta_N^{-i_1 h_1+\ell  i_1}\bigg(\eta_p(-1)+\eta_p(-a)\bigg)\\
&-M^2\sum_{h_2\in I}\sum_{i_2=0}^{N-1}\zeta_N^{-i_2 h_2+\ell i_2}\sum_{x\in \F_{q}\setminus\{-1\}}\chi_N^{i_2}(x+1)
+M^2\sum_{h_2\in I}\sum_{i_2\in A}\zeta_N^{-i_2 h_2+\ell i_2}(\eta_p(1)+\eta_p(-a+1))\\
&-M^2\sum_{h_3\in I}\sum_{i_3=0}^{N-1}\zeta_N^{-i_3 h_3+\ell i_3}\sum_{x\in \F_{q}\setminus\{-a\}}\chi_N^{i_3}(x+a)
+M^2\sum_{h_3\in I}\sum_{i_3\in A}\zeta_N^{-i_3 h_3+\ell  i_3}(\eta_p(a)+\eta_p(-1+a))\\
&+\sum_{h_1,h_2,h_3\in I}\sum_{i_1,i_2,i_3\in A}\zeta_N^{-i_1 h_1-i_2 h_2-i_3 h_3+\ell(i_1+i_2+i_3)}\sum_{x\in \F_{q}}\chi_N^{i_1}(x)\chi_N^{i_2}(x+1)\chi_N^{i_3}(x+a)+(q-3)M^3,
\end{align*}}
where $A=\{2j+1\,|\,0\le j\le (q-3)/2\}$. Then, by $\sum_{x\in \F_q^\ast}\chi_N(x)=0$ and 
\[
|D\cap (D+a)|=\frac{1}{N^2}\sum_{h,h'\in I}\sum_{i,i'=0}^{N-1}\zeta_N^{-i h-i' h'}\sum_{x\in \F_{q}\setminus\{0,a\}}\chi_N^{i}(x)\chi_N^{i'}(x-a), 
\] 
the above is also reformulated as 
\begin{align} 
&MN^2\left(|(D-\omega^\ell)\cap (D-a\omega^\ell)|+|D\cap (D-a\omega^\ell)|+
|D\cap (D-\omega^\ell)|\right)\nonumber\\
&-M(\eta_p(a)+\eta_p(-a+1)+\eta_p(a^2-a))\bigg(\sum_{h\in I}\sum_{i\in A}\zeta_N^{i(\ell-h)}\bigg)^2
-3M^3(q-1)\nonumber\\
&-M^2\left(\eta_p(-1)+\eta_p(-a)+\eta_p(1)+\eta_p(-a+1)+\eta_p(a)+\eta_p(-1+a)\right)\bigg(\sum_{h\in I}\sum_{i\in A}\zeta_N^{-i h+\ell i}\bigg)\nonumber\\
&+\sum_{h_1,h_2,h_3\in I}\sum_{i_1,i_2,i_3\in A}\zeta_N^{-i_1 h_1-i_2 h_2-i_3 h_3+\ell(i_1+i_2+i_3)}\sum_{x\in \F_{q}}\chi_N^{i_1}(x)\chi_N^{i_2}(x+1)\chi_N^{i_3}(x+a)+(q-3)M^3.\label{tripledef}
\end{align}
Let
\[
S_{i_1,i_2,i_3}(\omega^\ell,I)=\sum_{h_1,h_2,h_3\in I}\zeta_N^{-i_1 h_1-i_2 h_2-i_3 h_3+\ell(i_1+i_2+i_3)}\sum_{x\in \F_{q}}\chi_N^{i_1}(x)\chi_N^{i_2}(x+1)\chi_N^{i_3}(x+a).
\]
Recall that $D$ is a skew Hadamard difference set, i.e., $|D\cap (D+x)|=\frac{q-3}{4}$ for all $x\in \F_q^\ast$. 
Then, by $\eta_p(-1)=-1$, the sum (\ref{tripledef}) is reformulated as 
\begin{align*}
&\sum_{i_1,i_2,i_3\in A}S_{i_1,i_2,i_3}(\omega^\ell,I)+(q-3)M^3+3MN^2\frac{q-3}{4}
-3M^3(q-1)\\
&\hspace{1cm}-M(\eta_p(a)+\eta_p(-a+1)+\eta_p(a^2-a))\bigg(\sum_{h\in I}\sum_{i\in A}\zeta_N^{i(\ell- h)}\bigg)^2.
\end{align*}
Since $I\,(\mod{p_1^m})=\Z/p_1^m\Z$, there is exactly one 
$h'\in I$ such that $\ell-h'\equiv 0\,(\mod{p_1^m})$. Write $\ell-h'=\epsilon p_1^m$. 
Then, by noting that $A=\{2j+1\,|\,0\le j\le p_1^m-1\}$, we have 
\begin{align*}
\bigg(\sum_{h\in I}\sum_{i\in A}\zeta_N^{i(\ell- h)}\bigg)^2
=\bigg(\sum_{h\in I}\zeta_{N}^{\ell-h}\sum_{j=0}^{p_1^m-1}\zeta_{p_1^m}^{j(\ell- h)}\bigg)^2
=\bigg(p_1^m\zeta_{N}^{\epsilon p_1^m}\bigg)^2=p_1^{2m}.
\end{align*}
Thus, we finally have 
\begin{align}
N^3\cdot  T_{\omega^\ell,a}(D)&=\sum_{i_1,i_2,i_3\in A}S_{i_1,i_2,i_3}(\omega^\ell,I)+(q-3)M^3+3MN^2\frac{q-3}{4}
-3M^3(q-1)\nonumber\\
&\hspace{1cm}-M(\eta_p(a)+\eta_p(-a+1)+\eta_p(a^2-a))p_1^{2m}. \label{Tvalue}
\end{align}
\subsection{Triple intersection numbers modulo a prime extension degree}
The following theorem gives a sufficient condition for  lifts of  Feng-Xiang 
skew Hadamard difference sets being inequivalent to the Paley difference sets. 
\begin{theorem}\label{mainineq1}
Let $t$ be an odd prime with $\gcd{(t,p_1)}=1$. 
Let $D=\bigcup_{i \in I}C_{i}^{(N,q)}$ be a Feng-Xiang skew Hadamard difference set and $D'$ be the lift of $D$ to $\F_{q^t}$.  
If the set $\{T_{\omega^\ell,a}(D^{(t^{-1})})\,(\mod{t})\,|\,0\le \ell \le N-1\}$
contains $u$ distinct numbers, then so does 
$\{T_{\gamma^\ell,a}(D')\,(\mod{t})\,|\,0\le \ell \le N-1\}$, where $\omega$ and $\gamma$
are primitive elements of $\F_q$ and $\F_{q^t}$, respectively. 
\end{theorem}
\proof
Without loss of generality, we can assume that $\omega=\gamma^{(q^t-1)/(q-1)}$. 
Let $\chi_N$ be a multiplicative character of order $N$ of $\F_{q}$ such that $\chi_N(\omega)=\zeta_N$ and 
$\chi_N'$ be the lift of $\chi_N$ to $\F_{q^t}$. 
Define 
\[
S_{i_1,i_2,i_3}^{(t)}(\gamma^\ell,I)=\sum_{h_1,h_2,h_3\in I}\zeta_N^{-i_1 h_1-i_2 h_2-i_3 h_3+\ell(i_1+i_2+i_3)}\sum_{x\in \F_{q^t}}{\chi'}_N^{i_1}(x){\chi'}_N^{i_2}(x+1){\chi'}_N^{i_3}(x+a).
\]
Then, by Eq.~(\ref{Tvalue}), we have 
\begin{align*}
N^3\cdot  T_{\gamma^\ell,a}(D')=&\sum_{i_1,i_2,i_3\in A}S_{i_1,i_2,i_3}^{(t)}(\gamma^\ell,I)+(q^t-3)M^3+3MN^2\frac{q^t-3}{4}
-3M^3(q^t-1)\nonumber\\
&\hspace{1cm}-M(\eta_p(a)+\eta_p(-a+1)+\eta_p(a^2-a))p_1^{2m}. 
\end{align*}
By Eq.~(\ref{weil1}) of Theorem~\ref{multiweil}, there are two complex numbers $w_1,w_2$ such that 
\[
\sum_{x\in \F_{q}}{\chi}_N^{i_1}(x){\chi}_N^{i_2}(x+1){\chi}_N^{i_3}(x+a)=
-w_1-w_2
\]
and 
\[
\sum_{x\in \F_{q^t}}{\chi'}_N^{i_1}(x){\chi'}_N^{i_2}(x+1){\chi'}_N^{i_3}(x+a)=
-w_1^t-w_2^t. 
\]
Since $t$ is an odd prime satisfying $\gcd{(t,p_1)}=1$, we have 
\begin{align*}
\sum_{x\in \F_{q^t}}{\chi'}_N^{i_1}(x){\chi'}_N^{i_2}(x+1){\chi'}_N^{i_3}(x+a)
&=-w_1^t-w_2^t\\
&\equiv (-w_1-w_2)^t\,(\mod{t})\\
&\equiv \bigg(\sum_{x\in \F_{q}}\chi_N^{i_1}(x)\chi_N^{i_2}(x+1)\chi_N^{i_3}(x+a)\bigg)^t\,(\mod{t})\\
&\equiv \sum_{x\in \F_{q}}\chi_N^{t i_1}(x)\chi_N^{t i_2}(x+1)\chi_N^{t i_3}(x+a)\,(\mod{t})
\end{align*}
Therefore, we obtain 
\begin{align}
&\hspace{0.5cm} S_{i_1,i_2,i_3}^{(t)}(\gamma^\ell,I)\nonumber\\
&\equiv\sum_{h_1,h_2,h_3\in I}\zeta_N^{-ti_1(t^{-1} h_1)-ti_2(t^{-1} h_2)-ti_3( t^{-1}h_3)+(\ell\cdot t^{-1})(ti_1+ti_2+ti_3)}\sum_{x\in \F_{q}}\chi_N^{ti_1}(x)\chi_N^{ti_2}(x+1)\chi_N^{ti_3}(x+a)\,(\mod{t})\nonumber\\ 
&=S_{ti_1,ti_2,ti_3}(\omega^{t^{-1}\ell}, t^{-1}I). \nonumber
\end{align}
Hence,  it holds that  
\begin{align}
\bigg\{\sum_{i_1,i_2,i_3\in A}S_{i_1,i_2,i_3}^{(t)}(\gamma^\ell,I)\,|\,0\le \ell\le N-1\bigg\}&\equiv 
\bigg\{\sum_{i_1,i_2,i_3\in A}S_{ti_1,ti_2,ti_3}(\omega^{t^{-1}\ell},t^{-1}I)\,|\,0\le \ell\le N-1\bigg\}\, (\mod{t})\nonumber\\
&\equiv \bigg\{\sum_{i_1,i_2,i_3\in A}S_{i_1,i_2,i_3}(\omega^\ell,t^{-1}I)\,|\,0\le \ell\le N-1\bigg\}\,(\mod{t}).\label{maintheq2}
\end{align}
By the assumption that  the set $\{T_{\omega^\ell,a}(D^{(t^{-1})})\, (\mod{t})\,|\,0\le \ell \le N-1\}$ contains $u$ distinct numbers, 
we have  
\[
\bigg|\bigg\{\sum_{i_1,i_2,i_3\in A}S_{i_1,i_2,i_3}(\omega^\ell,t^{-1} I)\,\, \,  (\mod{t})\,|\,0\le \ell\le N-1\bigg\}\bigg|= u. 
\]
Then, by Eq.~(\ref{maintheq2}), we also have 
\[
\bigg|\bigg\{\sum_{i_1,i_2,i_3\in A}S_{i_1,i_2,i_3}^{(t)}(\gamma^\ell,I)\,\, \,   (\mod{t})\,|\,0\le \ell\le N-1\bigg\}\bigg|= u, 
\]
i.e., the set 
$\{T_{\gamma^\ell,a}(D')\,(\mod{t})\,|\,0\le \ell \le N-1\}$ contains $u$  distinct numbers. 
\qed

\begin{remark}\label{larget}{\em 
\begin{itemize}
\item[(i)] 
Let $D=\bigcup_{i \in I}C_{i}^{(N,q)}$ be defined as in Theorem~\ref{mainineq1}.   
Assume that $\{T_{\omega^\ell,a}(D^{(t^{-1})})\,|\,0\le \ell\le N-1\}$ contains $u(\ge 3)$ distinct numbers, say, $a_1<a_2<\cdots<a_u$.
Put \[
v=\min\{a_{j+2}-a_{j}\,|\,1\le j\le u-2\}. 
\]
Let $t$ be any odd prime satisfying $t>v$ and $\gcd{(t,p_1)}=1$  
and let $D'$ be the lift of $D$ to $\F_{q^t}$.  
Then, we have 
\[
\bigg|\{T_{\omega^\ell,a}(D^{(t^{-1})})\,(\mod{t})\,|\,0\le \ell\le N-1\}\bigg|\ge 3. 
\] 
Hence,  by Theorem~\ref{mainineq1}, we have 
$|\{T_{\gamma^\ell,a}(D')\,(\mod{t})\,|\,0\le \ell \le N-1\}|\ge 3$.   (More roughly, one may take $t$ so as $t>a_u-a_1$.) 
\item[(ii)] 
If two Feng-Xiang skew Hadamard difference sets $D_1$ and $D_2$ in $\F_q$ satisfy \[
\{T_{\omega^\ell,a}(D_1^{(t^{-1})})\,(\mod{t})\,|\,0\le \ell\le N-1\}\not=\{T_{\omega^\ell,a}(D_2^{(t^{-1})})\,(\mod{t})\,|\,0\le \ell\le N-1\},
\]
then by the proof of Theorem~\ref{mainineq1} 
their lifts $D_1'$ and $D_2'$ to $\F_{q^t}$ also satisfy 
\[
\{T_{\gamma^\ell,a}(D_1')\,(\mod{t})\,|\,0\le \ell\le N-1\}\not=\{T_{\gamma^\ell,a}(D_2')\,(\mod{t})\,|\,0\le \ell\le N-1\}, 
\]
i.e., $D_1'$ and $D_2'$ are inequivalent. 
\end{itemize}
}
\end{remark}
\begin{corollary}\label{Corrr}
Let $D=\bigcup_{i \in I}C_{i}^{(N,q)}$ be a Feng-Xiang skew Hadamard difference set. 
Let $t$ be any odd prime greater than $4N^3\sqrt{q}$   
and $D'$ be the lift of $D$ to $\F_{q^t}$.  
If $\{T_{\omega^\ell,a}(D^{(t^{-1})})\,|\,0\le i\le N-1\}$ contains $u$ distinct numbers, then so does 
$\{T_{\gamma^\ell,a}(D')\,(\mod{t})\,|\,0\le \ell \le N-1\}$, where $\omega$ and $\gamma$ are primitive elements of $\F_q$ and $\F_{q^t}$, respectively. 
\end{corollary}
\proof
Assume that $\{T_{\omega^\ell,a}(D^{(t^{-1})})\,|\,0\le i\le N-1\}$ contains $u$  distinct numbers. 
By Remark~\ref{larget}~(i), it is enough to show that $a_u-a_1\le 4N^3\sqrt{q}$. 
By Eq.~(\ref{Tvalue}), it is clear that 
\[
a_u-a_1\le\frac{2}{N^3}
\max\bigg\{\bigg|\sum_{i_1,i_2,i_3\in A}S_{i_1,i_2,i_3}(\omega^\ell,t^{-1}I)\bigg|: 0\le \ell\le N-1\bigg\}. 
\]
Now we estimate $|\sum_{i_1,i_2,i_3\in A}S_{i_1,i_2,i_3}(\omega^\ell ,t^{-1}I)|$. By Eq.~(\ref{weil2}) of Theorem~\ref{multiweil}, we have
\begin{align*}
\bigg|\sum_{i_1,i_2,i_3\in A}S_{i_1,i_2,i_3}(\omega^\ell,t^{-1}I)\bigg|&\le \sum_{i_1,i_2,i_3\in A}
\sum_{h_1,h_2,h_3\in t^{-1}I}\zeta_N^{-i_1 h_1-i_2 h_2-i_3 h_3+\ell(i_1+i_2+i_3)}\bigg|\sum_{x\in \F_{q}}\chi_N^{i_1}(x)\chi_N^{i_2}(x+1)\chi_N^{i_3}(x+a)\bigg|\\
&\le |A|^3|I|^3\sqrt{q}=2N^{6}\sqrt{q}. 
\end{align*}
Thus, we obtain 
\begin{align*}
a_u-a_1\le\frac{2}{N^3}
\max\bigg\{\bigg|\sum_{i_1,i_2,i_3\in A}S_{i_1,i_2,i_3}(\omega^\ell,t^{-1}I)\bigg|: 0\le \ell\le N-1\bigg\}\le 4 N^3 \sqrt{q}. 
\end{align*}
This completes the proof of the corollary. 
\qed

\vspace{0.5cm}
Corollary~\ref{Corrr} implies that for any sufficiently large prime $t$ the lift of  a Feng-Xiang skew Hadamard difference set $D$ in $\F_q$ to $\F_{q^t}$ is inequivalent to the Paley difference set if 
$|\{T_{\omega^\ell,a}(D^{(t^{-1})})\,|\,0\le \ell \le N-1\}|\ge 3$. 
\begin{example}\label{example1}{\em 
Let $p=11$, $N=2p_1=14$, $f=3$, and $I=\langle p\rangle \cup -2\langle p\rangle\cup \{0\}\,(\mod{N})$. Then, we have 
$I(\mod{p_1})=\Z/{p_1}\Z$ and the conditions of Theorem~\ref{FXtheorem} are  
satisfied. Thus, $D=\bigcup_{i\in I}C_{i}^{(N,p^f)}=\bigcup_{i\in I}\omega^i\langle \omega^N\rangle$ is a Feng-Xiang skew Hadamard difference set, where $\omega$ is a primitive element of $\F_{p^f}$. 
Now, we consider the triple intersection numbers $T_{\omega^\ell,a}(D^{(t^{-1})})$, $0\le \ell \le N-1$, with $a=3$. By Magma,  the author  checked that 
\[
\{T_{\omega^\ell,a}(D^{(t^{-1})})\,|\,0\le \ell \le N-1\}=\{147, 158, 164, 167, 173, 184\}
\]
for any odd $1\le t<N$ with  $\gcd{(t,p_1)}=1$. 
This implies that $D^{(t^{-1})}$ is inequivalent to the Paley difference set. 

It is clear that for any odd prime $t$ with  $\gcd{(t,p_1)}=1$ it holds that 
\[
\bigg|\{T_{\omega^\ell,a}(D^{(t^{-1})})\,(\mod{t})\,|\,0\le \ell \le N-1\}\bigg|\ge 3.  
\]
Hence, by Theorem~\ref{mainineq1}, the lift $D'$ of $D$ to $\F_{p^{ft}}$ for any odd prime $t$ with $\gcd{(p_1,t)}=1$ satisfies 
\[
\bigg|\{T_{\gamma^\ell,a}(D')\,(\mod{t})\,|\,0\le \ell \le N-1\}\bigg|\ge 3,
\] 
where $\gamma$
is a primitive element of $\F_{p^{ft}}$. Thus, $D'$ is also inequivalent to the Paley difference set. Furthermore, by applying Theorem~\ref{mainineq1} recursively, 
 the lift $ D''$ of $D$ to $\F_{p^{ft^{h}}}$ for any  $h\ge 1$ is also inequivalent to the Paley difference set.  The extension degrees $t^h$ less that $50$ covered by Theorem~\ref{mainineq1} in this case  are listed below: 
\[
t^h=3,5,7,9,11,13,17,19,23,25,27,29,31,37,41,43,47,49. 
\] 
}\end{example}
\subsection{Triple intersection numbers modulo the characteristic $p$}
The extension degree $t$ in Theorem~\ref{mainineq1} is limited to 
a prime. The following theorem allows us to take $t$ as an arbitrary odd integer. 
\begin{theorem}\label{maintheo2}
Let $t$ be any odd positive integer and consider the $p$-adic expansion $t=\sum_{h=0}^{r}x_hp^h$ with $0\le x_h\le p-1$. 
Write $e(t)=\sum_{h=0}^r x_h$ and then there is an odd integer $t'$ such that 
$t'=e(e(\cdots e(t)\cdots))$ and $1\le t'\le p-2$. Let $D=\bigcup_{i \in I}C_{i}^{(N,q)}$ be a Feng-Xiang skew Hadamard difference set and let $D'$ and 
$D''$ be its lifts to $\F_{q^t}$ and $\F_{q^{t'}}$, respectively. 
If the set $\{T_{\beta^\ell,a}(D'')\,(\mod{p})\,|\,0\le \ell \le N-1\}$
contains $u$ distinct numbers, then so does 
$\{T_{\gamma^\ell,a}(D')\,(\mod{p})\,|\,0\le \ell \le N-1\}$, where $\beta$ and $\gamma$ are primitive elements of $\F_{q^{t'}}$ and $\F_{q^t}$, respectively. 
\end{theorem}
\proof
Let $\omega$ be a primitive element of $\F_q$. Without loss of generality, we can assume that  $\omega=\beta^{(q^{t'}-1)/(q-1)}=\gamma^{(q^{t}-1)/(q-1)}$. 
Let $\chi_N$ be a multiplicative character of order $N$ of $\F_{q}$ such that $\chi_N(\omega)=\zeta_N$ and let 
$\chi_N'$ and $\chi_{N}''$ be the lifts of $\chi_N$ to $\F_{q^t}$ and $\F_{q^{t'}}$, respectively. 
Define  \[
S_{i_1,i_2,i_3}^{(t)}(\gamma^\ell,I)=\sum_{h_1,h_2,h_3\in I}\zeta_N^{-i_1 h_1-i_2 h_2-i_3 h_3+\ell(i_1+i_2+i_3)}\sum_{x\in \F_{q^t}}{\chi'}_N^{i_1}(x){\chi'}_N^{i_2}(x+1){\chi'}_N^{i_3}(x+a).
\]
Then, by Eq.~(\ref{Tvalue}), we have 
\begin{align*}
N^3\cdot  T_{\gamma^\ell,a}(D')=&\sum_{i_1,i_2,i_3\in A}S_{i_1,i_2,i_3}^{(t)}(\gamma^\ell,I)+(q^t-3)M^3+3MN^2\frac{q^t-3}{4}
-3M^3(q^t-1)\nonumber\\
&\hspace{1cm}-M(\eta_p(a)+\eta_p(-a+1)+\eta_p(a^2-a))p_1^{2m}. 
\end{align*}
By Eq.~(\ref{weil1}) of Theorem~\ref{multiweil}, there are two complex numbers $w_1,w_2$ such that 
\[
\sum_{x\in \F_{q}}{\chi}_N^{i_1}(x){\chi}_N^{i_2}(x+1){\chi}_N^{i_3}(x+a)=
-w_1-w_2, 
\]
\[
\sum_{x\in \F_{q^{t}}}{\chi'}_N^{i_1}(x){\chi'}_N^{i_2}(x+1){\chi'}_N^{i_3}(x+a)=
-w_1^t-w_2^t,
\]
and 
\[
\sum_{x\in \F_{q^{t'}}}{\chi''}_N^{i_1}(x){\chi''}_N^{i_2}(x+1){\chi''}_N^{i_3}(x+a)=
-w_1^{t'}-w_2^{t'}. 
\]
By Proposition~\ref{p-reduc}, we have 
\begin{align*}
\sum_{x\in \F_{q^t}}{\chi'}_N^{i_1}(x){\chi'}_N^{i_2}(x+1){\chi'}_N^{i_3}(x+a)
&\equiv \bigg(\sum_{x\in \F_{q}}\chi_N^{i_1}(x)\chi_N^{i_2}(x+1)\chi_N^{i_3}(x+a)\bigg)^t\,(\mod{p})\\
&= \prod_{h=0}^r\bigg(\sum_{x\in \F_{q}}\chi_N^{i_1}(x)\chi_N^{i_2}(x+1)\chi_N^{i_3}(x+a)\bigg)^{x_hp^h}\\
&\equiv \prod_{h=0}^r\bigg(\sum_{x\in \F_{q}}\chi_N^{p^h i_1}(x)\chi_N^{p^h i_2}(x+1)\chi_N^{p^h i_3}(x+a)\bigg)^{x_h}\,(\mod{p})\\
&=\prod_{h=0}^r\bigg(\sum_{x\in \F_{q}}\chi_N^{i_1}(x^{p^h})\chi_N^{i_2}(x^{p^h}+1)\chi_N^{i_3}(x^{p^h}+a)\bigg)^{x_h}\\
&=\prod_{h=0}^r\bigg(\sum_{x\in \F_{q}}\chi_N^{i_1}(x)\chi_N^{i_2}(x+1)\chi_N^{i_3}(x+a)\bigg)^{x_h}\\
&=\bigg(\sum_{x\in \F_{q}}\chi_N^{i_1}(x)\chi_N^{i_2}(x+1)\chi_N^{i_3}(x+a)\bigg)^{e(t)}. 
\end{align*}
By repeating this computation, we have 
\begin{align*}
\sum_{x\in \F_{q^t}}{\chi'}_N^{i_1}(x){\chi'}_N^{i_2}(x+1){\chi'}_N^{i_3}(x+a)&\equiv
\bigg(\sum_{x\in \F_{q}}\chi_N^{i_1}(x)\chi_N^{i_2}(x+1)\chi_N^{i_3}(x+a)\bigg)^{t'}\,(\mod{p})\\
&\equiv
\sum_{x\in \F_{q^{t'}}}{\chi''}_N^{i_1}(x){\chi''}_N^{i_2}(x+1){\chi''}_N^{i_3}(x+a)\, \, (\mod{p}). 
\end{align*}
Therefore, we obtain 
\begin{align*}
S_{i_1,i_2,i_3}^{(t)}(\gamma^\ell,I)&=\sum_{h_1,h_2,h_3\in I}\zeta_N^{-i_1 h_1-i_2 h_2-i_3 h_3+\ell(i_1+i_2+i_3)}\sum_{x\in \F_{q^t}}{\chi'}_N^{i_1}(x){\chi'}_N^{i_2}(x+1){\chi'}_N^{i_3}(x+a)\\
&\equiv\sum_{h_1,h_2,h_3\in I}\zeta_N^{-i_1 h_1-i_2 h_2-i_3 h_3+\ell(i_1+i_2+i_3)}\sum_{x\in \F_{q^{t'}}}{\chi''}_N^{i_1}(x){\chi''}_N^{i_2}(x+1){\chi''}_N^{i_3}(x+a)\, \, (\mod{p})\\
&=S_{i_1,i_2,i_3}^{(t')}(\beta^\ell,I)
\end{align*}
and 
\begin{equation}\label{maintheq3}
\bigg\{\sum_{i_1,i_2,i_3\in A}S_{i_1,i_2,i_3}^{(t)}(\gamma^\ell,I)\,|\,0\le \ell\le N-1\bigg\}\equiv 
\bigg\{\sum_{i_1,i_2,i_3\in A}S_{i_1,i_2,i_3}^{(t')}(\beta^\ell,I)\,|\,0\le \ell\le N-1\bigg\}\, (\mod{p}).
\end{equation}
By the assumption that  the set $\{T_{\beta^\ell,a}(D'')\, (\mod{p})\,|\,0\le \ell \le N-1\}$ contains $u$ distinct numbers, 
we have  
\[
\bigg|\bigg\{\sum_{i_1,i_2,i_3\in A}S_{i_1,i_2,i_3}^{(t')}(\beta^\ell,I)\, (\mod{p})\,|\,0\le \ell\le N-1\bigg\}\bigg|=u. 
\]
Then, by Eq.~(\ref{maintheq3}), we also have 
\[
\bigg|\bigg\{\sum_{i_1,i_2,i_3\in A}S_{i_1,i_2,i_3}^{(t)}(\gamma^\ell,I)\, (\mod{p})\,|\,0\le \ell\le N-1\bigg\}\bigg|=u,  
\]
i.e., the set 
$\{T_{\gamma^\ell,a}(D')\,(\mod{t})\,|\,0\le \ell \le N-1\}$ contains $u$ distinct numbers. 
\qed

\vspace{0.5cm}
Theorem~\ref{maintheo2} implies that for any odd integer $t$ the lift $D'$ of a Feng-Xiang skew Hadamard difference set $D$ in $\F_q$ to $\F_{q^t}$ is inequivalent to the Paley difference set if 
the lift $D''$ of $D$ to $\F_{q^{t'}}$ satisfies 
$|\{T_{\beta^\ell,a}(D'')\,(\mod{p})\,|\,0\le \ell \le N-1\}|\ge 3$ for every odd $1\le t'\le p-2$.  
\begin{example}{\em 
Let $p$, $N$, $f$, $a$, and $I$ be defined as in Example~\ref{example1}. Then $D=\bigcup_{i\in I}C_i^{(N,q)}$  satisfies 
\[
\{T_{\omega^\ell,a}(D)\,|\,0\le \ell \le N-1\}=\{147, 158, 164, 167, 173, 184\}. 
\] 
It is clear that 
\[
\bigg|\{T_{\omega^\ell,a}(D)\,(\mod{p})\,|\,0\le \ell \le N-1\}\bigg|\ge 3.  
\]
Hence, by Theorem~\ref{maintheo2}, the lift $D'$ of $D$ to $\F_{p^{ft}}$ for any odd $t$ with $e(e(\cdots e(t)\cdots))=1$ satisfies 
\[
\bigg|\{T_{\gamma^\ell,a}(D')\,(\mod{p})\,|\,0\le \ell \le N-1\}\bigg|\ge 3,
\] 
where $\gamma$
is a primitive element of $\F_{p^{ft}}$. Thus, $D'$ is inequivalent to the Paley difference set. The extension degrees $t>1$ less than $50$ covered by Theorem~\ref{maintheo2} in this case are $
t=11,21,31,$ and $41$. Note that $t=21$ is not covered by Theorem~\ref{mainineq1}. 
}\end{example}
\section{Concluding remarks}\label{Sec4}
In this paper, we obtained two theorems which give sufficient conditions for 
lifts of  Feng-Xiang 
skew Hadamard difference sets being inequivalent to the Paley difference sets. 
As an example, we showed that there are infinitely many integers $t$ such that 
the lifts of the Feng-Xiang difference set $D=\bigcup_{i\in \langle 11\rangle \cup -2\langle 11\rangle\cup\{0\}}C_i^{(14,11^3)}$ to $\F_{11^{3t}}$ are inequivalent to the Paley difference sets. Further small examples are listed  in Table~\ref{Tab1}. (In the table, let $\omega$ be a primitive element of $\F_{p^f}$ and $n_t:=|\{T_{\omega^\ell,3}(D^{(t^{-1})})\,(\mod{t})\,|\,0\le \ell \le N-1\}|$.) In these examples, we fixed parameters as $N=14$ and $f=3$ due to the memory-capacity of computer.  For each $p\in\{11,23,67,79,107\}$, the Feng-Xiang skew Hadamard difference sets $D=\bigcup_{i\in I}C_i^{(N,p^f)}$ for $I$ listed in the table and their lifts to $\F_{p^{ft}}$ for sufficiently large odd primes $t$ are inequivalent to the Paley difference sets. Moreover, the lifts are mutually inequivalent  by  Remark~\ref{larget}~(ii). 
{\tiny 
\begin{table}[t]
\caption{Examples of Feng-Xiang
skew Hadamard difference sets 
and their triple intersection numbers
}
\label{Tab1}
$$
\begin{array}{|c||c|c|c|}
\hline
\mbox{$(p,f,N)$}&\mbox{index set $I$}&\mbox{$\{T_{\omega^\ell,3}(D^{(t^{-1})})\,|\,0\le \ell \le N-1\}$ for  $\gcd{(t,N)}=1$}&n_t\mbox{ for  odd primes $t$} \\
\hline \hline
&\{ 0, 1, 2, 3, 4, 5, 6 \}&\{ 159, 162, 164, 167, 169, 172 \}&\mbox{$n_5=2$ and $n_t\ge 3$ for any other $ t$} \\ \cline{2-4}
(11,3,14) &\{0, 1, 2, 3, 4, 6, 12 \}&\{ 157, 160, 165, 166, 171, 174 \}&\mbox{$n_3=2$ and $n_t\ge 3$ for any other  $ t $}\\ \cline{2-4}
 &\{0,1, 6,9,10,11, 12\}&\{ 147, 158, 164, 167, 173, 184 \}&\mbox{$n_t\ge 3$ for any $ t$}\\ \cline{2-4}
& \{0, 1, 2, 4, 6, 10, 12\}&\{ 163, 164, 167, 168 \}&\mbox{$n_t\ge 3$ for any $ t\ge 3 $}\\ \hline
(11,3,2) &\{0\}\, \, \mbox{(Paley)}&\{157, 174\}&\mbox{$n_{17}=1$ and $n_t= 2$ for any other $ t $}\\
\hline \hline 
&\{ 0, 1, 2, 3, 4, 5, 6 \}&\{ 1497, 1498, 1503, 1515, 1525, 1537, 1542, 1543 \}&\mbox{$n_3=2$ and $n_t\ge 3$ for any other $ t $}\\
\cline{2-4}
&\{0, 1, 2, 3, 4, 6, 12 \}&\{  1498, 1503, 1508, 1514, 1526, 1532, 1537, 1542 \}&\mbox{$n_t\ge 3$ for any $ t $}\\ \cline{2-4}
(23,3,14)  &\{0,1, 6,9,10,11, 12\}&\{ 1481, 1509, 1514, 1526, 1531, 1559  \}&\mbox{$n_5=2$ and $n_t\ge 3$ for any other  $ t $}\\ \cline{2-4}
& \{0, 1, 2, 4, 6, 10, 12\}&\{1508, 1514, 1526, 1532 \}&\mbox{$n_3=1$ and  $n_t\ge 3$ for any $ t $}\\\hline
(23,3,2) &\{0\}\, \, \mbox{(Paley)}&\{1520\}&\mbox{$n_t=1$ for any $ t $}\\
\hline \hline
&\{ 0, 1, 2, 3, 4, 5, 6 \}&\hspace{-1cm}\{ 37457, 37519, 37525, 37587,&\mbox{$n_t\ge 3$ for any $ t $}\\
&&\hspace{1cm} 37602, 37664, 37670, 37732\}&\\
 \cline{2-4}
 &\{0, 1, 2, 3, 4, 6, 12 \}&\hspace{-1cm}\{37453, 37523, 37587, 37591, &\mbox{$n_t\ge 3$ for any $ t $}\\
(67,3,14)  &&\hspace{1cm} 37598, 37602, 37666, 37736  \}&\\
\cline{2-4}
&\{0,1, 6,9,10,11, 12\}&\{ 37526, 37587, 37594, 37595, 37602, 37663  \}&\mbox{$n_t\ge 3$ for any $ t $}\\
 \cline{2-4}
& \{0, 1, 2, 4, 6, 10, 12\}&\{37543, 37559, 37630, 37646 \}&\mbox{$n_3,n_{29}=2$ and $n_t\ge 3$ for any other $ t $}\\ \hline
(67,3,2) &\{0\}\, \, \mbox{(Paley)}&\{37502, 37687\}&\mbox{$n_t= 2$ for any $ t $}\\
\hline \hline
&\{ 0, 1, 2, 3, 4, 5, 6 \}&\hspace{-1cm}\{ 61470, 61575, 61607, 61623, &\mbox{$n_t\ge 3$ for any $ t $}\\
& &\hspace{1cm} 61636, 61652, 61684, 61789 \}&\\
 \cline{2-4}
 &\{0, 1, 2, 3, 4, 6, 12 \}&\hspace{-1cm}\{  61398, 61535, 61549, 61552, &\mbox{$n_t\ge 3$ for any $ t $}\\
(79,3,14)  &&\hspace{1cm} 61707, 61710, 61724, 61861\}&\\
 \cline{2-4}
&\{0,1, 6,9,10,11, 12\}&\{61513, 61533, 61546, 61713, 61726, 61746 \}&\mbox{$n_3,n_{5}=2$ and $n_t\ge 3$ for any other $ t $}\\
 \cline{2-4}
& \{0, 1, 2, 4, 6, 10, 12\}&\{61434, 61511, 61748, 61825 \}&\mbox{$n_7,n_{11},n_{157}=2$ and $n_t\ge 3$ for any other $ t $}\\
\hline
(79,3,2) &\{0\}\, \, \mbox{(Paley)}&\{ 61519, 61740 \}&\mbox{$n_t= 2$ for any $ t $}\\
\hline \hline
&\{ 0, 1, 2, 3, 4, 5, 6 \}&\hspace{-1cm}\{152751, 152895, 152976, 153021,&\mbox{$n_3=2$ and $n_t\ge 3$ for any other $ t $}\\
& &\hspace{1cm} 153238, 153283, 153364, 153508  \}&\\ \cline{2-4}
(107,3,14)  &\{0, 1, 2, 3, 4, 6, 12 \}&\{ 152969, 153065, 153092, 153167, 153194, 153290\}&\mbox{$n_3=1$ and $n_t\ge 3$ for any other $ t $}\\ \cline{2-4}
 &\{0,1, 6,9,10,11, 12\}&\{152643, 153040, 153102, 153157, 153219, 153616  \} &\mbox{$n_3=2$ and $n_t\ge 3$ for any other $ t $}\\
\cline{2-4}
& \{0, 1, 2, 4, 6, 10, 12\}&\{ 153028, 153103, 153156, 153231\}&\mbox{$n_3,n_5=2$ and $n_t\ge 3$ for any other $ t $}\\
\hline
(107,3,2) &\{0\}\, \, \mbox{(Paley)}&\{152977, 153282\}&\mbox{$n_t=2$ for any $ t $}\\
\hline
\end{array}
$$
\end{table}}

Our main theorems works well for Feng-Xiang skew Hadamard difference sets since the difference sets have the nice property that their lifts are
again skew Hadamard difference sets. Here, we have the following natural question: 
are there skew Hadamard difference sets with the ``lifting'' property 
other than Paley difference sets and Feng-Xiang skew Hadamard difference sets? 
Below, we give an immediate generalization of Feng-Xiang skew Hadamard difference sets. 
\begin{theorem}\label{GeneFX}
Let $p_1$ be a prime, $N=2p_1^m$, and 
let $p\equiv 3\,(\mod{4})$ be a prime such that $2\in \langle p\rangle\, (\mod{p_1^m})$, $\gcd{(p_1,p-1)}=1$, and $f:=\ord_N(p)$ is odd. 
Let 
$s$ be any odd integer, $I$ any subset of $\Z/N\Z$ such that 
$\{i\,(\mod{p_1^m})\,|\,i\in I\}=\Z/p_1^m\Z$, and  let $D=\bigcup_{i\in I}C_i^{(N,q)}$, where $q=p^{fs}$. Then, $D$ is a skew Hadamard 
difference set. 
\end{theorem}
The above theorem follows immediately from \cite[Theorem 11.3.5]{BEW97} and the proof of \cite[Theorem 3.2]{FX113}. This can be seen as follows. 
By Theorem 11.3.5 (the Davenport-Hasse product formula on Gauss sums) in \cite{BEW97} and the assumption that $2\in \langle p\rangle\, (\mod{p_1^m})$, we have 
\[
G_{f}(\chi_N)=\frac{G_f(\chi_{p_1^m})G_f(\chi_2)}{\chi_{p_1^m}(2)G_f(\chi_{p_1^m}^{2^{-1}})}=G_f(\chi_2)\chi_{p_1^m}^{-1}(2). 
\]
By $\gcd{(p_1,p-1)}=1$, the restriction of $\chi_{p_1^m}$ to $\F_p$ is 
trivial. Hence, we have $G_{f}(\chi_N)=G_f(\chi_2)$. The remaining proof is 
same with that of Theorem 3.2 in \cite{FX113}. 
For example, the case where $(p,p_1,f)=(47,127,21)$ is covered by Theorem~\ref{GeneFX} but not covered by Theorem~\ref{FXtheorem}. 

Interesting problems which are worth looking into as future works are listed below. 
\begin{enumerate}
\item[(1)] We showed that there are infinitely many skew Hadamard difference sets inequivalent to the Paley difference sets by using  ``recursive'' techniques not saying
anything about the inequivalence of  ``starting'' skew Hadamard difference sets theoretically.  Determine whether the ``starting'' skew Hadamard difference sets are inequivalent to the Paley difference sets without using computer.  
\item[(2)] Find skew Hadamard difference sets with the lifting property, i.e., their lifts are again skew Hadamard difference sets,  
other than the Paley difference sets and the difference sets of Theorem~\ref{GeneFX}. 
\item[(3)] Recently, several new constructions of skew Hadamard difference sets have been known other than Feng-Xiang skew Hadamard difference sets \cite{CF,DAW,DY06,DWX07,Momi,WH09}. Determine whether such skew Hadamard difference sets are inequivalent to the Paley difference sets.  
\end{enumerate}

\end{document}